\documentclass[11pt,leqno]{article}

\usepackage{amssymb,amsmath,amsthm}

\theoremstyle{plain}
\newtheorem{thm}{Theorem}
\newtheorem{cor}[thm]{Corollary}

\theoremstyle{remark}

\theoremstyle{definition}

\newcommand{\vp}{\varepsilon}
\newcommand{\ph}{\varphi}

\newcommand{\n}{\noindent}

\def\CC{ \;Ê{}^{ {}_\vert }\!\!\!{\rm C}}

\def\NN{{\rm I}\!{\rm N}}

\theoremstyle{remark}
\newtheorem{rem}{Remark}
\begin{document}

\title{A Remark on Hypercontractive Semigroups and Operator Ideals}

\author{by\\ \phantom{blank line}\\ Gilles Pisier\footnote{Partially supported by  
NSF grant 
No.~0503688}\\\\ Department of Mathematics\\ 
Texas A\&M University\\ College Station, TX 77843-3368, 
U. S. A.
\\and\\ Universit\'e Paris VI\\ Equipe d'Analyse, Case 186, 75252\\ 
Paris Cedex 05, France\\  {\bf email: pisier@math.tamu.edu}}

\date{August 2007}
\maketitle

\abstract{ In this note, we answer a question raised by Johnson and Schechtman \cite{JS}, about the hypercontractive semigroup on $\{-1,1\}^{\NN}$. More generally, we prove the theorem below.  }
\bigskip

\n 2000 Math. Subject Classification: Primary 47D20.
Seconday 47L20.

\bigskip

In this note, we answer a question raised by Johnson and Schechtman \cite{JS}. We refer to their paper for motivation and background. In addition, we refer the reader to \cite{Hyp,Hyp1,Hyp2,Hyp3} for  examples of hypercontractive semigroups and related background (e.g. the classical works
 of Bonami, Nelson,  Federbush, Gross, Beckner,...).

\begin{thm} Let $1<p<2$.
Let $(T(t))_{t>0}$ be a holomorphic semigroup
on $L_p$ (relative to a probability space).
Assume the following mild
form of hypercontractivity: for some  large enough number $s>0$, 
$T(s)$ is bounded from $L_p$ to $L_2$. Then for any $t>0$, $T(t)$ is in the norm closure  in $B(L_p)$ 
(denoted by $\overline{\Gamma_2}$)  of the subset 
(denoted by ${\Gamma_2}$) formed by the operators mapping $L_p$ to $L_2$ (a fortiori these operators factor through a Hilbert space).
More precisely, for any $t>0$ there are constant $c,d$ such that for any $\vp>0$ there is an operator $T':\ L_p \to L_p$ with $\|T(t)-T'\|\le \vp$ such that $\|T':\ L_p \to L_2\|\le c\vp^d.$
\end{thm}
\begin{proof} Unless specified otherwise, by the norm
of $T(t)$ we mean its norm as acting from $L_p$ to itself.
 We will consider 
 the compact set $K$ such that its boundary is a triangle  with vertices
 $0$, $a\pm ib$, where $a,b>0$ are chosen 
 (by holomorphy of $T$) so that 
$t\mapsto T(t)$ is bounded on $K$. By the semigroup property,  $t\mapsto T(t)$ is also bounded
on $s+K$. Moreover, for any $t\in K$, 
the factorization $T(s+t)=T(t)T(s)$ shows that
$T(s+t)\in  {\Gamma_2}$. More precisely, by assumption
we have constants $c,d$ such that
$$\|T(s):\ L_p\to L_2\|= c<\infty$$  and 
$$\sup_{t\in K}\|T(t):\ L_p\to L_p\|= d<\infty.$$ Therefore,
 there is a constant $C_1=cd$ such that
$$\forall t\in K\quad \|T(s+t):\ L_p\to L_2\|\le C_1.$$
Consider now the open domain $V$ bounded by the triangle with vertices $0$ and $s+a\pm ib$. The boundary of $V$ consists of two parts: the vertical segment joining $s+a\pm ib$ that we denote by $V_1$ and the remaining part of the boundary that we denote by $V_0$. Note that, on one hand,  $t\mapsto T(t)$ is bounded on $V_0$ in $B(L_p)$; more precisely,
if $t=w+a+i w' \in V_0$, then, since
$  T(t) = T(w) T(a+i w')$, we have   $$\sup_{t\in V_0}\|  T(t) \|\le c\sup_{0<w<s}\|T(w)\|=C_0.\leqno(1)$$ On the other hand, 
 we have
$$\sup_{t\in V_1}\|T(t):\ L_p\to L_2\|\le d\|T(s):\ L_p\to L_2\|=C_1.\leqno(2)$$ To complete the proof it suffices to show that $T(t)\in   \overline{\Gamma_2}$ for any $0<t<s$. We will show this more generally for any $t$ in the interior of $V$.
For any such $t$, let $\mu$ be the harmonic measure of $V$ with respect to the interior point $t$.
Recall that $\mu$ is a probability on the boundary of $V$ such that for any bounded analytic function $f$ on $V$, we have
$$f(t)=\int_{\partial V} f(z) \mu(dz).\leqno(3)$$
We will now use classical ideas from complex interpolation ({\it cf.} \cite{BeL}), modulo the conformal equivalence of $V$ with the standard  vertical strip in the complex plane. We may rewrite $\mu$ as
$\mu= (1-\theta)  \mu_0+  \theta \mu_1$ where
$\mu_0$ and $\mu_1$ are probabilities on respectively $V_0$ and $V_1$.

Let $S=\{ z\in \CC\mid 0<\Re(z)<1$ and
let $S_j=\{ z\in \CC\mid  \Re(z)=j$ ($j=0,1$).

\n {\bf Claim 1:} By the Riemann mapping theorem, there is a conformal mapping from $V$ onto $S$
taking $t$ to $\theta$, and taking the inner part of $V_j$ to $S_j$ ($j=0,1$). Indeed, $V$ is conformally equivalent to the unit disc $D$
by a mapping $\ph$  
taking $t$ to $0$ 
 that is also a homeomorphism from $\bar V$ to $\bar D$ (see {\it e.g.} \cite[p.114]{Ra}). Therfore, $\ph$  must take  $\mu$ to the normalized Lebesgue measure on the unit circle (that is
the harmonic measure of $D$ with respect  to
the point $\ph(t)=0$). Since $\ph$ takes the boundary of $V$ continuously  to the unit circle,  it takes $V_0$ (resp.  $V_1$)  to an arc
$I_0$ (resp.  $I_1$)  of the unit circle of measure $1-\theta$ (resp.  $\theta$). But now it is well known that
there is a conformal equivalence $\chi_\theta$ from $S$ to $D$
taking $S_j$ to $I_j$ ($j=0,1$) and mapping $\theta$
to $0$. 
Indeed, quite explicitly: the function $\chi_\theta\colon\ S\to D$ defined by
$$\chi_\theta(z)= \frac{e^{i\pi z} - e^{i\pi \theta}}{e^{i\pi z} - e^{-i\pi \theta}}$$ does the job (up to a rotation).
This proves our claim.

\n {\bf Claim 2:} For any $\vp>0$ there is a function
$\psi\in H^\infty(V)$ such that $\psi(t)=1$
and such that $\sup_{z\in V_0}|\psi(z)|\le \vp$, and
 $\sup_{z\in V_1}|\psi(z)|\le \vp^{(\theta-1)/\theta}$.
 
 Indeed, by conformal equivalence, it suffices to produce a function $\xi\in H^\infty(S)$
 such that $\xi(\theta)=1$
 and $\sup_{z\in S_0}|\xi(z)|\le \vp$, and
 $\sup_{z\in S_1}|\xi(z)|\le \vp^{(\theta-1)/\theta}$. But then the classical function 
 $\xi(z)= \vp^{(\theta-z)/\theta}$
 (related to  geometric means), has the desired property.

 We now come to the conclusion. Fix $\vp>0$.
 By (3) we can write   $$T(t)=\psi(t)T(t)=\int \psi(z)T(z)  \mu(dz)$$ and hence
 $$ T(t)=(1-\theta)T_0+ \theta T_1$$
 where
 $$T_0=\int_{V_0} \psi (z) T(z) \mu_0(dz),\quad 
 T_1=\int_{V_1} \psi (z)T(z) \mu_1(dz).$$
 Finally, by  (1), (2) and Claim 2,    we have
 $$\|T_0\|\le C_0\vp\quad \|T_1:\ L_p\to L_2\|\le C_1 \vp^{(\theta-1)/\theta}.$$
 Thus we obtain $T_1\in \Gamma_2$ such that
 $\|T(t)-T_1\|\le C_0\vp$. Since $\vp$ is arbitrary, this implies $T(t)\in \overline{\Gamma_2}$.
\end{proof}
The next statement 
extends the corresponding statement in \cite{JS}. Of course, this is obvious if $X$ is such that   $T(t)_{|X}$ is diagonal, the point is that it is valid for an arbitrary subspace.

\begin{cor} In the above situation, if $T(t)$ restricted to a subspace $X\subset L_p$ is an isomorphism, then
$X$ must be isomorphic to a Hilbert space (and $X$ is complemented in $L_p$).
\end{cor}

\begin{proof} By perturbation, we may assume that 
there is a $T$ in ${\Gamma_2}$ such that  $T_{|X}$ is an isomorphism. This clearly  implies that $X$ is isomorphic to a Hilbert space. Moreover,   the orthogonal projection from $L_2$ onto $T(X)\subset L_2$ composed with $T^{-1}:\ T(X)\to X$ is a bounded projection from $L_p$ onto $X$.
\end{proof}
{\bf Remarks:} The above proof uses  a very general  (simple) principle, the same argument applies to more general situations e.g. to the hypercontractive Fermionic semigroup of Carlen-Lieb \cite{CL}.
Moreover, one can also consider hypercontractivity
on $H^p$ following Svante Janson's and Weissler's results \cite{Ja, We}.
Then the main result is valid for $0<p<2$. In that case, the semigroup is formed of compact operators, so, in the situation of the corollary, the subspace $X$ must be finite dimensional.

More generally, the same proof as above yields: 
\begin{thm}
Let $L$ be  any Banach space, and let
$T(t)$ be a holomorphic semigroup in $B(L)$. Let
$\Gamma\subset B(L)$ be a left ideal equipped
with a Banach space norm $\gamma$ such that, for some constant $C$, we have
$$\forall x\in B(L)\ \ \forall T\in \Gamma\quad \gamma(xT)\le C\|x\|_{B(L)} \gamma(T).$$
Assume that $T(s)$ is in $\Gamma$ for some $s>0$.
Then, $T(t)$ is in $\overline{\Gamma}$, the norm closure 
of  $\Gamma$ in $B(L)$  for all $t>0$. 
More precisely, for any $t>0$ there are constant $c,d$ such that for any $\vp>0$ there is an operator $T'\in B(L) $ with $\|T(t)-T'\|\le \vp$ such that $\gamma(T')\|\le c\vp^d.$

\begin{rem} Wolfgang Arendt kindly pointed out to me that the holomorphy of the semigroup itself is enough for the first part of the conclusion of either Theorems above and that  this part is absolutely immediate: Indeed, say with the notation of the preceding statement, consider the quotient map
$Q:\ B(L)\to B(L)/\overline{\Gamma}$, then $t\mapsto QT(t)$ vanishes for large $t$, and hence, by analyticity, must vanish identically. 
\end{rem}

\end{thm}

\end{document}